%% file: main.tex
\documentclass{amsart}
\input epsf

\usepackage{fullpage}

\usepackage{amsfonts,amsthm,amsmath,amssymb,latexsym}
\usepackage{graphicx,color}
\usepackage[all]{xy}

\def\D{{\bf D}}

\def\S{\mathcal{H}}

\begin{document}

\author{Sylvain Bonnot, R. C. Penner and Dragomir \v Sari\' c}

\address{Institute for Mathematical Sciences, Stony Brook University,
Stony Brook, NY 11794-3660} \email{bonnot@math.sunysb.edu}

\address{Departments of Mathematics and Physics/Astronomy,
University of Southern California, Los Angeles, CA 90089}
\email{rpenner@math.usc.edu}

\address{Institute for Mathematical Sciences, Stony Brook University,
Stony Brook, NY 11794-3660} \email{saric@math.sunysb.edu}

\title[Presentation for the mapping class group of the solenoid]
{A presentation for the baseleaf preserving mapping class group of
the punctured solenoid}

\subjclass{}

\keywords{}

\begin{abstract}
We give a presentation for the baseleaf preserving mapping class
group $Mod(\S )$ of the punctured solenoid $\S$. The generators
for our presentation were introduced previously, and several relations
among them were derived.  In addition, we
show that $Mod(\S )$ has no non-trivial central elements. Our main
tool is a new complex of triangulations of the disk upon which $Mod(\S )$ acts.

\end{abstract}

\maketitle

\thispagestyle{empty}
\input{imsmark}
\SBIMSMark{2006/4}{July 2006}{}

\section{Introduction}
This note continues the investigation (begun in \cite{PS}) of the
baseleaf preserving mapping class group $Mod(\S )$ for the
punctured solenoid $\S$. Our main result is a presentation for
$Mod(\S )$.  The punctured solenoid $\S$ is an inverse limit of
the system of all finite unbranched covers of a punctured surface
of negative Euler characteristic, and its baseleaf preserving
mapping class group $Mod(\S )$ consists of all homotopy classes of
appropriate self maps of $\S$ which preserve a distinguished leaf.  (See Section
2 or \cite{PS} for more details.) One motivation for studying $Mod(\S
)$ comes from Sullivan's observation \cite{NS} that the Ehrenpreis
conjecture is equivalent to the statement that $Mod(\S )$ has
dense orbits in the Teichm\"uller space of the solenoid $\S$. Another motivation
is that the baseleaf preserving mapping class group
$Mod(\S )$ is a large subgroup of the studied commensurator group
$Comm(F_2)$ of the free group $F_2$ on two generators (for
the definition, see Section 2). Namely, if we identify $F_2$ with the
once punctured torus group $G$, then $Mod(\S )$ is the subgroup of
$Comm(G)$ which preserves peripheral elements, i.e., preserves parabolic elements.

\vskip .2 cm

The path components of the punctured solenoid $\S$ are called {\it
leaves} and the {\it baseleaf} is a fixed distinguished leaf.
Since $Mod(\S )$ preserves the baseleaf, which is dense
in $\S$, it is enough to analyse the action on the baseleaf. We
recall that the baseleaf (and indeed any other leaf) is conformally
equivalent to the unit disk $\D$.

\vskip .2cm

Given an ideal triangulation of
the unit disk (i.e., the baseleaf) which is invariant under a
finite-index subgroup $K$ of $PSL_2(\mathbb{Z})$ and a specified
edge of the triangulation, there are two adjacent triangles which
together form a  ``neighboring'' quadrilateral.  We may replace the specified edge
of this quadrilateral by its other diagonal, and performing  this modification
for each edge in the $K$-orbit of the specified edge, we define the $K$-equivariant
Whitehead move.
The resulting ideal triangulation is also invariant
under $K$. A Whitehead homeomorphism of $S^1$ is obtained by
mapping an ideal triangulation of the unit disk onto its image
under a Whitehead move. It is shown in \cite{PS} that the
Whitehead homeomorphisms together with $PSL_2(\mathbb{Z})$
generate the baseleaf preserving mapping class group $Mod(\S )$
(see Section 2 or \cite{PS} for more details). In fact, $Mod(\S)$
consists of quasisymmetric homeomorphisms of $S^1$ which
conjugate one finite-index subgroup of $PSL_2(\mathbb{Z})$ onto
another \cite{PS} or \cite{Od}. In \cite{PS}, four relations among
Whitehead homeomorphisms are identified, and three of them arise
in our presentation (see Theorem 4.4(c).)

\vskip .2 cm

We first introduce the {\it triangulation complex} $\mathcal X$
for the punctured solenoid $\S$. The vertices of $\mathcal X$ are
TLC tesselations, i.e., ideal triangulations of ${\bf D}$ invariant under some finite-index subgroup of $PSL_2(\mathbb{Z})$. Two
vertices of $\mathcal X$ are joined by an edge if they differ by a
Whitehead move.

\vskip .2cm

There are several types of two-cells in $\mathcal X$: two edges of a triangulation
may have disjoint neighboring quadrilaterals, in which case there is a two-cell corresponding
to commutativity of their associated Whitehead moves; the two edges may have neighboring
quadrilaterals which share a triangle, in which case there is a two-cell corresponding to
the pentagon relation; or a single Whitehead move equivariant for a finite-index subgroup $K<PSL_2({\mathbb Z})$
may be written as the finite composition of Whitehead moves equivariant for a subgroup
of $K$ of finite index.  (The two-cells are described more precisely in Section 3).

\vskip .2 cm

\paragraph{\bf Theorem 3.1} {\it The triangulation complex $\mathcal X$
is connected and simply connected.}

\vskip .2 cm

The action of $Mod(\S )$ on the triangulation complex $\mathcal X$
is evidently cellular.
Furthermore (see \cite{PS}), there is only one orbit of vertices in $\mathcal X$,
and the isotropy group of a vertex $v$, i.e., its stabilizer $\Gamma (v)$, is a conjugate of  $PSL_2(\mathbb{Z})$.
Together with the further analysis of the isotropy group $\Gamma (E)$ of an unoriented edge $E$, standard techniques \cite{Bro}
allow us to derive a presentation
of  $Mod(\S )$. To simplify this
presentation, we actually choose a larger set of generators for $Mod(\S )$, namely,
we take as generators all Whitehead moves starting from the basepoint of
$\mathcal X$.
(This is a smaller set of generators than in  \cite{PS} but larger than necessary.)
We
denote by $\mathcal{E}^{+}$ the set of edges of $\mathcal X$ which contain the
basepoint and are not inverted by an element of $Mod(\S )$, and
by $\mathcal{E}^{-}$ the set of edges which contain the basepoint
and are inverted by an element of $Mod(\S )$.
It is necessary to
fix one Whitehead homeomorphism $g_E$ for each edge
$E\in\mathcal{E}^{+}$ in a consistent way.  (See Section 4 regarding
this choice.)  Let $\mathcal{E}^{\pm}= {\mathcal E}^+\sqcup {\mathcal E}^-$
denote the set of unoriented edges, and let $\Gamma ^+(E)$ denote the subgroup of
$\Gamma (E)$ which does not invert the edge $E\in{\mathcal E}^-$.

\vskip .2 cm

\paragraph{\bf Theorem 4.3} {\it The modular group $Mod(\S )$ is
generated by the isotropy subgroup $PSL_2(\mathbb{Z})$ of the
basepoint $\tau_{*}\in\mathcal X$, the isotropy subgroups
$\Gamma (E)$ for $E\in\mathcal{E}^{\pm}$,
and by the elements $g_E$ for
$E\in\mathcal{E}^{+}$. The following relations on these
generators give a complete presentation of $Mod(\S )$:

\vskip .1 cm

\noindent a) The inclusions of $\Gamma (E)$ into $PSL_2(\mathbb{Z})$,
for $E\in\mathcal{E}^{+}$, are given by $\Gamma (E)=K'$,
where the terminal endpoint of $E$ is invariant under
the finite-index subgroup $K'<PSL_2(\mathbb{Z})$;

\vskip .1 cm

\noindent b) The inclusions of $\Gamma ^+(E)$ into $PSL_2(\mathbb{Z})$,
for $E\in\mathcal{E}^{-}$, are given by $\Gamma (E)=K'$,
where the terminal endpoint of $E$ is invariant under
the finite-index subgroup $K'<PSL_2(\mathbb{Z})$;

\vskip .1 cm

\noindent c) The relations introduced by the boundary edge-paths
of two-cells in $\mathcal F$ given by the equations}
(\ref{pentagon}), (\ref{pentagon1}), (\ref{pentagon2}),
(\ref{square}), (\ref{square1}), (\ref{coset}) {\it and}
(\ref{coset1});

\vskip .1 cm

{\it \noindent d) The redundancy relations: for any two edges $E$
and $E'$ in $\mathcal{E}^{\pm}$ and for any $\gamma\in
PSL_2(\mathbb{Z})$ such that $\gamma (E)=E'$, we get the relation
$$g_{E'}\circ \gamma '=\gamma\circ g_E,$$
 where $\gamma '$ is the unique
element of $PSL_2(\mathbb{Z})$ that satisfies $\gamma '(e_0)=e_1'$
with $e_1'=g_{E'}^{-1}( \gamma (e_0))$.}

\vskip .2 cm

It is well-known that the mapping class group of a Riemann surface of finite type
has trivial center provided the genus is at least three, and we obtain the analogous result
for $Mod(\S )$.

\vskip .2 cm

\paragraph{\bf Theorem 5.1} {\it The modular group $Mod(\S )$ of
the punctured solenoid $\S$ has trivial center.}

\vskip .2 cm

Define $\mathcal{Y}=\mathcal{X}/Mod(\S )$ and let $\mathcal N$
be the subgroup of $Mod(\S )$ generated by all elements which fix
a point in $\mathcal X$.  By a standard result \cite{Arm}, we get

\vskip .2 cm

\paragraph{\bf Theorem 5.4} {\it The topological fundamental group of
${\mathcal Y}={\mathcal X}/Mod(\S )$ is given by
$$
\pi_1({\mathcal Y})= Mod(\S )/\mathcal{N}.
$$}

\vskip .2 cm

\paragraph{\it Acknowledgements} We are grateful to John Milnor
for useful comments.

\section{Preliminaries}

Fix a punctured surface $S$ ({\it the base surface}) with
negative Euler characteristic and empty boundary, and consider the system of all
finite unbranched covers of $S$. There is a partial
ordering on the covers as follows. If one cover $\pi_1$ can
be factored as the composition of two covers $\pi _1=\pi\circ\pi _2$, where $\pi, \pi _2$ are
also finite unbranched covers, then
$\pi_1\geq \pi_2$. The system of covers
is inverse directed, and there is thus an inverse limit.

\vskip .2 cm

\paragraph{\bf Definition 2.1} The {\it punctured solenoid} $\S$
is the inverse limit of the system of finite unbranched covers of
a punctured surface without boundary and with negative Euler characteristic.

\vskip .2 cm

The inverse limit does not depend on the base surface as long as
it is of negative Euler characteristic \cite{Od}, \cite{PS}. The
punctured solenoid $\S$ is locally homeomorphic to a disk times a
Cantor set. Each path component is called a {\it leaf} , and each leaf
is homeomorphic to the unit disk. The punctured solenoid $\S$ has
uncountably many leaves, each of which  is dense in $\S$. If we
require in the above definition of $\S$ that each punctured
surfaces and each covering map is pointed, we obtain a
distinguished point, called the {\it basepoint} of $\S$. The leaf
containing the basepoint is called the {\it baseleaf}. The
punctured solenoid $\S$ is a non-compact topological space with
one end, which is homeomorphic to a horoball
times a Cantor set modulo the continuous action of a countable group.
For more details, see
\cite{PS}.

\vskip .2 cm

\paragraph{\bf Definition 2.2} The {\it baseleaf preserving
mapping class group} $Mod(\S )$ of the punctured solenoid $\S$ is
the group of isotopy classes all self homeomorphisms of $\S$ which preserve the
baseleaf and which are quasiconformal on leaves.

\vskip .2 cm

The restriction of an element of $Mod(\S )$ to the baseleaf gives
a quasiconformal homeomorphism of the unit disk $\D$ (upon fixing an
identification of the baseleaf with $\D$) up to isotopy. Thus, an
element of $Mod(\S )$ determines a well-defined quasisymmetric homeomorphism of
$S^1$, and we shall thus identify $Mod(\S)$ with a appropriate
group of quasisymmetric maps.  (See Theorem~2.8.)

\vskip .2 cm

\paragraph{\bf Definition 2.3} The {\it commensurator group}
$Comm(G )$ of a group $G $ consists of
equivalence classes of isomorphisms of finite-index subgroups of
$G$, where two isomorphisms are equivalent if they agree on a finite-index subgroup in the intersection of their domains.

\vskip .2 cm

\paragraph{\bf Theorem 2.4} \cite{PS} {\it The modular group $Mod(\S )$ is
isomorphic to a proper subgroup of the commensurator group
$Comm(F_2)$ of the free group $F_2$ on two generators. Namely,
$Mod(\S )$ is isomorphic to the subgroup of $Comm(F_2)$ consisting of all elements
which preserve the peripheral elements under some fixed
identification $F_2\equiv G$, where $G<PSL_2(\mathbb{Z})$ is the
group uniformizing the once-punctured torus.}

\vskip .2 cm

In fact, it is convenient in the definition of the punctured solenoid $\S$
to fix the base surface to be the once-punctured torus $\D /G$,
where $G<PSL_2(\mathbb{Z})$. Given an isomorphism of two finite-index subgroups $K,H$ of $G$
which preserves peripheral elements,
there exists a unique quasisymmetric map of $S^1$ which conjugates
$K$ onto $H$. Thus, by the previous theorem, we may consider $Mod(\S )$
as a group of quasisymmetric maps of $S^1$ which conjugate one
finite-index subgroup of $PSL_2(\mathbb{Z})$ onto another.

\vskip .2 cm

We recall that the decorated Teichm\"uller space $\tilde{T} (\S )$
of the punctured solenoid $\S$ is partitioned into sets
according to the bending information of the
convex hull construction \cite{PS}.  When the bending locus is a
triangulation on each leaf, then it is locally constant in the
transverse direction, and the action of the baseleaf preserving
mapping class group $Mod(\S )$ is transitive on this subspace of $\tilde{T}(\S)$.
It is convenient to consider the action of $Mod(\S )$ on the ideal
triangulations of the baseleaf (i.e., the unit disk $\D$) arising
by restrictions (to the baseleaf) of the triangulations of $\S$.

\vskip .2 cm

\paragraph{\bf Definition 2.5} A {\it transversely locally
constant (TLC) tesselation} $\tau$ of the unit disk $\D$ is
a lift to $\D$ of an
ideal triangulation of some punctured surface $\D/K$ of finite type, i.e.,
an
ideal triangulation of $\D$ invariant under a finite-index
subgroup $K$ of $PSL_2(\mathbb{Z})$, where the ideal points of the tesselation agree with
$\bar{\mathbb{Q}}\subset S^1$.

\vskip .2 cm

A particularly important example of a TLC tesselation is the Farey
tesselation $\tau_{*}$ (see, for example, \cite{P1} or \cite{PS}), which is invariant
under the group $PSL_2(\mathbb{Z})$. Let $K$ be a finite-index
subgroup of $PSL_2(\mathbb{Z})$ and let $\tau$ be a $K$-invariant
TLC tesselation of $\D$. A {\it characteristic map} for $\tau$ is
a homeomorphism $h:S^1\to S^1$ such that $h(\tau_{*})=\tau$ (see
\cite{P1}).

\vskip .2cm

The map $h=h(\tau ,e)$ is completely determined by specifying an oriented edge
$e\in\tau$, namely, the standard
oriented edge $e_0=(-1,1)$ in $\tau_{*}$ is mapped onto $e$,
the triangle to the left or right of $e_0$ in $\tau_*$ is mapped to the triangle to the left
or right, respectively, of $e$ in $\tau$, and so on.
Note that any two
characteristic maps for $\tau$ differ by pre-composition with an
element of $PSL_2(\mathbb{Z})$.

\vskip .2 cm

\paragraph{\bf Theorem 2.6} \cite{PS} {\it The characteristic map
$h=h(\tau ,e)$ for a $K$-invariant TLC tesselation $\tau$
conjugates a finite-index subgroup $H$ of $PSL_2(\mathbb{Z})$ onto
$K$.}

\vskip .2 cm

Note that a characteristic map necessarily conjugates peripherals
into peripherals because it is a homeomorphisms of $S^1$, and so
a characteristic map for a TLC tesselation lies in $Mod(\S )$.

\vskip .2 cm

Let $\tau$ be a $K$-invariant TLC tesselation of $\D$ with a
distinguished oriented edge $e$; we allow for the possibility that  $\tau$ is
invariant under a larger subgroup of $PSL_2(\mathbb{Z})$. Fix an
edge $f\in\tau$. We form a new $K$-invariant TLC tesselation $\tau
'$ by replacing each $\gamma (f)$, for $\gamma\in K$, by $\gamma
(f')$, where $f'$ is the diagonal of the unique ideal
quadrilateral in $(\D -\tau )\cup\{ f\}$ different from $f$. We
say that $\tau '$ is obtained from $\tau$ by performing a {\it
Whitehead move} along $K\{ f\}$. If $e\notin K\{ f\}$ then we let
$e$ be the distinguished oriented edge in $\tau '$ as well;
if $e=\gamma (f)$ for some $\gamma\in K$, then we let
$e'=\gamma (f')$ be the distinguished oriented edge for $\tau '$,
where $e'$ is given the orientation such that the tangent vectors
to $e$ and $e'$ at their intersection point comprise a positively oriented
basis for the oriented disk $\D$.

\vskip .2 cm

\paragraph{\bf Definition 2.7} Let $\tau$ be a $K$-invariant TLC
tesselation of $\D$ with a distinguished oriented edge $e$. The
{\it Whitehead homeomorphism} for $\tau$ and $e$ is
$$
k(\tau ,e)=h(\tau ',e')\circ h(\tau ,e)^{-1},
$$
where $(\tau ',e')$ arises from $(\tau,e)$ under the
Whitehead move and $h(\tau ',e'), h(\tau ,e)$ are the characteristic
maps.

\vskip .2 cm

A Whitehead homeomorphism lies in $Mod(\S
)$ since it is  the composition of two elements of
the group $Mod(\S )$ by Theorem~2.6.

\vskip .2 cm

\paragraph{\bf Theorem 2.8}\cite{PS} {\it The modular group $Mod(\S )$ of the
punctured solenoid $\S$ is generated by $PSL_2(\mathbb{Z})$ and by
Whitehead homeomorphisms for all TLC tesselations. In addition,
$Mod(\S )$ acts transitively on the set of all TLC tesselations of the
unit disk $\D$.}

\section{The triangulation complex}

We introduce a two-complex $\mathcal X$ associated to TLC
tesselations of $\D$. This complex $\mathcal X$
is an adaptation to our situation of the two-skeleton of the complex
dual to the cell decomposition of the decorated
Teichm\"uller space introduced by Penner \cite{P2} and
Harer \cite{Ha}. On the other hand, $\mathcal X$
is analogous to the complex of cut systems
of Hatcher and Thurston \cite{HT} in that
$Mod(\S )$ acts transitively on its vertices.

\vskip .2 cm

We begin the definition of $\mathcal X$ by giving its vertices.
A {\it vertex} of the {\it triangulation complex} ${\mathcal
X}$ is a TLC tesselations of the unit disk $\D$. The {\it
basepoint} of $\mathcal X$ is Farey tesselation $\tau_{*}$. A
characteristic map between any two TLC tesselations is an element
of $Mod(\S )$ \cite[Lemma 7.5]{PS}, and so $Mod(\S )$ acts
transitively on vertices of $\mathcal X$.

\vskip .2 cm

We next introduce edges of $\mathcal X$ with one endpoint at the
basepoint $\tau_{*}$. An unordered pair of vertices $\{\tau
,\tau_{*}\}$ determines an edge in $\mathcal X$ if $\tau$ can be
obtained from $\tau_{*}$ by a single Whitehead move, i.e., $\tau$
is obtained from $\tau_{*}$ by replacing an orbit $K\{ f\}$ of an
edge $f$ in $\tau_{*}$ by the orbit $K\{ f'\}$, where $K$ is
torsion free and of finite-index in $PSL_2(\mathbb{Z})$, and $f'$
is the diagonal of the quadrilateral in $(\D -\tau_{*})\cup\{ f\}$
different from $f$.

More generally and by definition, an unordered pair of vertices $\{\tau_1,\tau_2\}$
determines
an edge of $\mathcal X$ if $\{\tau_1,\tau_2\}$ is the image by an
element of $Mod(\S )$ of an edge $\{\tau_{*},\tau \}$ defined
above. In particular, this implies that if a TLC tesselation
$\tau_2$ is obtained by performing a $K$-invariant Whitehead move
on a TLC tesselation $\tau_1$ then $\{\tau_1,\tau_2\}$ is an edge
in $\mathcal X$. To see this, take a characteristic map $h$ for
$\tau_1$ (i.e., $h$ is a homeomorphism of $S^1$ such that
$h(\tau_{*})=\tau_1$) and consider $\{
\tau_{*},h^{-1}(\tau_2)\}$. Thus, $\tau_1$
and $\tau_2$ differ only in the $K$ orbit of the diagonals of an
ideal rectangle. The characteristic map $h$ conjugates a finite
index subgroup $H$ of $PSL_2(\mathbb{Z})$ onto $K$ \cite[Lemma
7.5]{PS}, so $h^{-1}(\tau_1)=\tau_{*}$ and $h^{-1}(\tau_2)$
differ only in that they have different diagonals on an $H$-orbit of a rectangle.
It follows that  $\{ \tau_{*},h^{-1}(\tau_2)\}$ is an edge
corresponding to an $H$-invariant Whitehead move, i.e.,
$\{\tau_1,\tau_2\} =h(\{
\tau_{*},h^{-1}(\tau_2)\})$ is an edge in $\mathcal
X$.

\vskip .2cm

However, there are edges which appear away from the basepoint
$\tau_{*}$ that do not correspond to Whitehead moves. They can be
described as a {\it generalized } Whitehead move invariant under a conjugate by a characteristic map of a finite-index subgroup of
$PSL_2(\mathbb{Z})$. Note that the set of edges in $\mathcal X$ is
invariant under the action of $Mod(\S )$ by construction.
This completes the definition of the one-skeleton of $\mathcal X$.

\vskip .2 cm

We introduce two-cells of $\mathcal X$ by first defining those that have
one vertex at the basepoint $\tau_{*}$.   Let $K$ be a torsion free finite-index
subgroup of $PSL_2(\mathbb{Z})$.  There are three types of two-cells:

\vskip .2cm

\noindent {\bf [Pentagon]}~ Suppose that $K$ is of index at least 9, i.e.,
$\tau_{*}/K$ is a triangulation of $\D /K$ which has at least
three complementary ideal triangles.   Any three adjacent
complementary triangles form a pentagon on $\D /K$ whose boundary
sides are possibly identified in pairs. Let $e_1$ and $e_2$ be two
representatives in $\tau_{*}$ of the diagonals of a
pentagon on $\D /K$ which share an ideal point. The sequence of five Whitehead
homeomorphisms $h(K,e_1)$, $h(K,e_2)$, $h(K,e_1')$, $h(K,e_2')$
and $h(K,e_1'')$ defines a closed edge-path in $\mathcal X$
based at $\tau_{*}$, where $e_1'$ is the new edge corresponding to
$e_1$ under the Whitehead move for $\tau_{*}$ along $K\{ e_1\}$,
$e_2'$ is the new edge corresponding to $e_2$ under the Whitehead
move for $\tau_1=h(K, e_1)(\tau_{*})$ along $K\{ e_2\}$, and
$e_1''$ is the new edge corresponding to $e_1'$ under the
Whitehead move for $\tau_2=h(K, e_2)(\tau_1)$ along $K\{ e_1'\}$
(see \cite{PS} for more details). We add a two-cell in
$\mathcal{X}$ whose boundary is this closed edge-path of length
five starting and ending at the basepoint $\tau_{*}$ and
call this two-cell a {\it pentagon} at the basepoint
$\tau_{*}$.

\vskip .2 cm

\noindent {\bf [Square]}~Let $K$ be a torsion free finite-index subgroup of
$PSL_2(\mathbb{Z})$ such that the triangulation $\tau_{*}/K$ of
$\D /K$ has two edges which do not lie in the boundary of a common complementary triangle.
Let $e_1$ and $e_2$
be two lifts to $\tau_{*}$ of the two non-adjacent edges.
Consider the closed edge-path of length four given by Whitehead
homeomorphisms $h(K,e_1)$, $h(K,e_2)$, $h(K,e_1')$ and
$h(K,e_2')$, where $e_1'$ corresponds to $e_1$ under the Whitehead
move $h(K,e_1)$ and $e_2'$ corresponds to $e_2$ under the
Whitehead move $h(K,e_2)$. We add a two-cell to $\mathcal X$ with
boundary equal to the above edge-path of length four and call it a
{\it square} cell at the basepoint $\tau_{*}$.

\vskip .2 cm

\noindent {\bf [Coset]}~Suppose $e\in\tau_{*}$ and let $H$ be a finite-index
subgroup of $K$. The orbit $K\{ e\}$ is canonically
decomposed into finitely many orbits $H\{ e_1\}, H\{ e_2\},\ldots
,H\{ e_k\}$, where $e_1=e, e_2,\dots ,e_k\in K\{ e\}$ and
$k=[K:H]$. Let $f$ be the other diagonal in the unique ideal
quadrilateral in $(\D -\tau_{*})\cup\{ e\}$ and let
$f_1,f_2,\ldots ,f_k\in K\{ f\}$ be the altered edges
corresponding to $e_1,e_2,\ldots ,e_k$. Consider a finite edge-path
based at $\tau_{*}$ consisting of the Whitehead
homeomorphisms $h(K,e)$, $h(H,f_1)$, $h(H,f_2)$, $\ldots$,
$h(H,f_k)$ corresponding to the tesselations $\tau_{*}$,
$\tau_1=h(K,e)(\tau_{*})$, $\tau_2=h(H,e_1)(\tau_1)$, $\ldots$,
$\tau_{k}=h(H,e_{k-1})(\tau_{k-1})$,
$\tau_{k+1}=h(H,e_k)=\tau_{*}$. We add a two-cell to $\mathcal X$
whose boundary is this edge-path and call it the {\it coset} cell
at the basepoint. Note that a different ordering of $f_1,\ldots
,f_k$ gives a different edge-path and hence a different coset
cell. In fact, there are $k!$ corresponding coset cells when
$[K:H]=k$. The edge $\{\tau_{*},\tau_1\}$ is called a {\it long
edge}, and all other edges are called {\it short edges}
corresponding to this coset cell.

\vskip .2 cm

Note that all two-cells introduced above have their boundaries
given by compositions of Whitehead moves invariant under subgroups of
$PSL_2(\mathbb{Z})$ as opposed to more general edges in $\mathcal
X$ where moves are only conjugate to Whitehead moves invariant
under subgroups of $PSL_2(\mathbb{Z})$.

\vskip .2 cm

To complete the definition of $\mathcal X$, an arbitrary two-cell in $\mathcal X$ is the image under
$Mod(\S )$ of a two-cell at the basepoint. If $h\in Mod(\S )$ and
$P$ is a two-cell based at $\tau_{*}$, then we say that $h(P)$
is based at $\tau :=h(\tau_{*})$. Note that closed edge-paths
based at $\tau_{*}$ are mapped to closed edge-paths, and hence
the boundaries of two-cells are well defined. The boundary of a pentagon or a
square two-cell
based at a tesselation $\tau\neq\tau_{*}$ each of whose edges is a
Whitehead move invariant with respect to a fixed finite-index
subgroup $K$ of $PSL_2(\mathbb{Z})$ (or equivalently, whose
vertices are TLC tesselations invariant under $K$) is likewise the boundary
of a two-cell.
Furthermore, the boundary of a coset two-cell starting at
$\tau\neq\tau_{*}$ whose initial vertex is invariant under $K$ and
whose other vertices are invariant under a subgroup $K_1<K$ of finite-index
is the image of the boundary of a  coset two-cell based at
$\tau_{*}$ by simply noting that a characteristic map which sends
$\tau_{*}$ onto $\tau$ conjugates $H_1<H$ onto
$K_1<K$ where $H<PSL_2(\mathbb{Z})$.

\vskip .2cm

By construction, the set of two-cells in $\mathcal X$
is invariant under $Mod(\S )$, and $Mod(\S )$ consequently acts
cellularly on the two-complex $\mathcal X$.

\vskip .2 cm

We claim that a pentagon or a square two-cell $P$ based at
$\tau\neq\tau_{*}$ with one vertex at the basepoint $\tau_{*}$ has
all edges given by Whitehead moves invariant under a fixed finite-index subgroup $K$ of $PSL_2(\mathbb{Z})$. Since any
characteristic map $h$ of $\tau_{*}$ onto the base vertex $\tau$
of the two-cell $P$ conjugates a finite-index subgroup
$H<PSL_2(\mathbb{Z})$ onto a finite-index subgroup
$K<PSL_2(\mathbb{Z})$ (under which $\tau$ is invariant), it follows that a two-cell (pentagon or square) $P'$ invariant under $H$ and based at
$\tau_{*}$ is mapped by $h$ onto the above two-cell $P$ based at
$\tau$ whose vertices are invariant under $K$. Moreover,
a coset two-cell $P$ which is based at $\tau\neq\tau_{*}$ whose
initial vertex is $\tau_{*}$ does not necessarily have edges arising
from Whitehead moves invariant under $K_1<K$. Let $P'$ be a coset
cell based at $\tau_{*}$ such that $h(P')=P$. In fact, if the
image under $h$ of the long edge of $P'$ is not incident on
$\tau_{*}$, then it is represented by a generalized Whitehead move
(invariant under a conjugate by $h$ of a subgroup of
$PSL_2(\mathbb{Z})$ which is not itself a group of M\"obius transformations).

\vskip.2 cm

\paragraph{\bf Theorem 3.1} {\it The triangulation complex $\mathcal X$
is connected and simply connected.}

\vskip .2 cm

\paragraph{\bf Proof} We first prove that $\mathcal X$ is
connected by showing that any vertex $\tau$ can be
connected to the basepoint $\tau_{*}$ by a finite edge-path. Let
$K$ be a finite-index subgroup of $PSL_2(\mathbb{Z})$ under which
$\tau$ is invariant. Thus, $\tau /K$ and $\tau_{*}/K$ are two
tesselations of a punctured surface $\D /K$. By results of
Penner \cite{P2} (or Harer \cite{Ha} or Hatcher \cite{Hat}), there
is a sequence of Whitehead moves on $\D /K$ which transforms
$\tau_{*}/K$ into $\tau /K$. The lifts of the Whitehead moves on
$\D /K$ to $\D$ are TLC Whitehead moves on $\D$ and they provide an
edge-path from $\tau_{*}$ to $\tau$ in $\mathcal X$. This establishes that
$\mathcal X$ is connected.

\vskip .2 cm

It remains to show that $\mathcal X$ is simply connected. We
recall a result of Harer \cite[Theorem 1.3]{Ha} or Penner \cite{P2}
for triangulations
of punctured surfaces. The set of top-dimensional simplices of the
triangulation complex of a finite punctured surface consists of
ideal triangulations, the codimension-one simplices are ideal
triangulations with one ideal geodesic erased, the codimension-two simplices are ideal triangulations with two ideal geodesics
erased, etc.. The main fact is that the triangulation complex of
the finite surface minus simplices which are given by
decompositions of the surface, where at least one complementary
component is not topologically a disc is homeomorphic to the
decorated Teichm\"uller space of the punctured surface. In particular, the
triangulation complex for a punctured surface is contractible.

\vskip .2 cm

Consider a closed edge-path $\alpha$ in the triangulation
complex $\mathcal X$ for the punctured solenoid. It is possible
that an edge $E$ in the path $\alpha$ is given by a generalized
Whitehead move, i.e., the two tesselations at the endpoints of $E$
are invariant under $hKh^{-1}$, where $K$ is a finite-index
subgroup of $PSL_2(\mathbb{Z})$ and $h\in Mod(\S )$. Since $h\in
Mod(\S )$, there exists $H_1,H_2<PSL_2(\mathbb{Z})$ of finite
index such that $H_2=hH_1h^{-1}$. Thus, $H_1\cap K=:K_1$ is of
finite-index in $K$, and we consider a coset two-cell corresponding to
the groups $K_1<K$ with long edge $h^{-1}(E)$. The edge
$h^{-1}(E)$ is homotopic modulo its endpoints to the path of short
edges in the coset two-cell, where each vertex is invariant under
$K_1$. The image under $h$ is a coset two-cell with long edge
corresponding to $hKh^{-1}$ and short edges corresponding to
$hK_1h^{-1}<PSL_2(\mathbb{Z})$. Thus, we can replace the long edge
invariant under $hKh^{-1}$, which is not a subgroup of
$PSL_2(\mathbb{Z})$, by the homotopic edge-path invariant under
$hK_1h^{-1}<H_2<PSL_2(\mathbb{Z})$. We may therefore replace $\alpha$
by an edge-path $\alpha '$ each of whose edges corresponds to a
Whitehead move invariant under a finite-index subgroup of
$PSL_2(\mathbb{Z})$ using only coset two-cells.

\vskip .2 cm

Let $K_1,K_2,\ldots ,K_n$ be finite-index subgroups of
$PSL_2(\mathbb{Z})$ which correspond to invariant Whitehead moves
defining the edges of $\alpha '$. Using coset two-cells
corresponding to each $K_i$, we first homotope the above edge-path
$\alpha '$ into a closed edge-path $\alpha ''$ where each edge
corresponds to an invariant Whitehead move with respect to a
single finite-index subgroup $K:=K_1\cap K_2\cap\ldots \cap K_n$.
The new edge-path $\alpha ''$ invariant under $K$ in the
triangulation complex $\mathcal X$ of the punctured solenoid can
be represented by a closed path $\gamma$ in the above
triangulation complex of a finite surface $\D /K$. (Recall that
$\mathcal X$ is an extension of the dual of the triangulation
complex of finite surface $\D /K$.) The path $\gamma$ starts and
ends in the top-dimensional simplex which corresponds to the
triangulation of $\D /K$ obtained by projecting the TLC
tesselation of $\D$ defining the basepoint of $\alpha''$ onto $\D
/K$. Furthermore, $\gamma$ crosses transversely codimension-one
simplices of the triangulation complex of $\D /K$ corresponding to
each edge in $\alpha ''$, and it enters each top-dimensional
simplex which corresponds to a vertex of $\alpha''$ in the given
order. Since the triangulation complex for punctured surface $\D
/K$ is simply connected \cite{Ha} or \cite{P2}, there exists a
homotopy of $\gamma$ into the trivial path which transversely
crosses codimension-two cells. The number of times the
homotopy crosses codimension-two cells is finite, and it is
possible to choose a homotopy which does not intersect
simplices of codimension greater than two. For each intersection point of the homotopy
with a codimension-two
simplex, there is a corresponding two-cell in $\mathcal X$ because
two-cells corresponding to ideal triangulations of the surface
with two edges erased lift to ideal triangulations of $\D$ with
orbits of two edges erased such that each complementary region is
finite sided. This exactly correspond to two-cells (either
pentagon or square) in $\mathcal X$. Thus, the homotopy for
$\gamma$ gives a homotopy between $\alpha''$ and the trivial path
in $\mathcal X$, and $\mathcal X$ is therefore simply connected. $\Box$

\vskip .2 cm

In the spirit of Ivanov's work \cite{Iva}, we may ask:

\vskip .2 cm

\paragraph{\bf Question} Is the group of automorphisms
$Aut(\mathcal{X})$ of the triangulation complex $\mathcal X$
isomorphic to the (extended) baseleaf preserving modular group?

\section{Presentation of $Mod(\S )$}

Applying a general theorem of Brown \cite[Theorem 1]{Bro} to the action of $Mod(\S)$ on $\mathcal X$,
we give a presentation for the modular group $Mod(\S )$ of the
punctured solenoid $\S$.

\vskip .2 cm

First recall \cite{PS} that $PSL_2(\mathbb{Z})$ is the isotropy group of the basepoint
$\tau_{*}\in\mathcal X$, that $Mod(\S)$ acts transitively on the vertices of
$\mathcal X$,
 that an arbitrary vertex
$\tau\in \mathcal X$ has isotropy group
$hPSL_2(\mathbb{Z})h^{-1}$, where the characteristic map $h:\tau_{*}\mapsto\tau$ lies in $Mod(\S)$,
and that $hPSL_2(\mathbb{Z})h^{-1}$ contains a
finite-index subgroup of $PSL_2(\mathbb{Z})$.

\vskip .2 cm

Consider the isotropy group of an edge in $\mathcal X$. Since
each vertex is mapped to the basepoint $\tau_{*}$, it is enough to
consider edges with an endpoint at $\tau_{*}$. The isotropy group
of any other edge is the conjugate of the isotropy group of such an
edge.

\vskip .2 cm

Let $E=\{\tau_{*},\tau\}$ be an arbitrary edge of $\mathcal X$
with one endpoint at the basepoint $\tau_{*}$ of $\mathcal X$.
There are two possibilities: either the isotropy group  $\Gamma (E)$
of $E$ contains elements which reverse the orientation of $E$
(i.e., interchanges $\tau_{*}$ and $\tau$), or
each element of $\Gamma (E)$ fixes each endpoint of $E$.

\vskip .2 cm

Let $\tau$ be obtained by a Whitehead move on $\tau_{*}$ invariant
under a torsion free finite-index subgroup $K$ of
$PSL_2(\mathbb{Z})$, and let us choose a characteristic map
$g:\tau_{*}\mapsto\tau$, where $g\in Mod(\S )$. Denote by $K'$ the
maximal extension of $K$ in $PSL_2(\mathbb{Z})$ which fixes
$\tau$. If $h\in \Gamma (E)$ preserves the orientation of
$E$ as above, then $h$ fixes both $\tau_{*}$ and $\tau$. By \cite[Lemma
7.3]{PS}, $h\in PSL_2(\mathbb{Z})$ and similarly $h\in
gPSL_2(\mathbb{Z})g^{-1}$, and so $h\in PSL_2(\mathbb{Z})\cap
gPSL_2(\mathbb{Z})g^{-1}=K'$. It follows that the subgroup
$\Gamma ^{+}(E)$ of the isotropy group $\Gamma (E)$ of an edge $E$ which
consists of elements which do not reverse orientation on $E$ is
equal to $K'<PSL_2(\mathbb{Z})$.

\vskip .2 cm

If $k\in \Gamma (E)$ reverses orientation of
$E=\{\tau_{*},\tau\}$, i.e., $k(\tau_{*})=\tau$ and $k(\tau
)=\tau_{*}$, then $k^2\in K'=\Gamma ^{+}(E)$. In particular, $k^2$ is a
lift of a self map of the Riemann surface $\D /K$. By
\cite[Proposition 1.3.6]{Odd}, $k$ is a lift of a self map of a
Riemann surface which finitely covers $\D /K$. We show that $k$ is
actually a lift of a self map of $\D /K$ itself.

\vskip .2 cm

\paragraph{\bf Lemma 4.1}
{\it Let $\tau$ be the image of $\tau_{*}$ under a $K$-invariant
Whitehead move, let
 $E=\{\tau_{*},\tau\}$ be the corresponding edge and
let $k\in \Gamma (E)-\Gamma ^{+}(E)$. Then $k$ conjugates $K$ onto
itself.}

\vskip .2 cm

\paragraph{\bf Proof} The proof proceeds in several steps.

\vskip .2 cm

\paragraph{\bf Simplification of the homeomorphism $k$.}
Since $k$ preserves the union of the two
tesselations $\tau, \tau_{*}$, it therefore sends a pair of
intersecting edges to a pair of intersecting edges.
Fix such an intersecting pair $e\in\tau_{*}$ and
$f\in \tau$, and consider the corresponding Whitehead move. Since $k(e)\in K\{
f\}$, there exists $\gamma\in K$ such that $(k\circ\gamma )(e)=f$.
It is enough to prove the lemma for $k\circ \gamma$ and we continue to denote it
by $k$.

\vskip .2 cm

\paragraph{\bf Orientation of edges}
Let us choose an orientation of $e$ and assign an orientation to
$f$ such that $k:e\mapsto f$ is orientation preserving. Assign
an orientation to each edge in the orbits of $e$ and $f$ under $K$
as follows. Let $e'=\gamma'(e)$ for some $\gamma'\in G$. Let
$\alpha$ be a differentiable arc connecting $e$ to $e'$ which
transversely crosses the minimal number of edges of $\tau_{*}$.
Give the orientation
to the curve $\alpha$ such that the tangent vector to $\alpha$ and
the tangent vector to $e$ at their point of  intersection form
a positively oriented basis for the tangent space of $\D$ at the
intersection point, and assign an orientation on $e'$ such that
the tangent vector to $\alpha$ and the tangent vector to $e'$ at
the intersection point $\alpha\cap e'$ form a positively oriented
basis to the tangent space. We may assign an orientation to any
$f'=\gamma' (f)$ in a similar fashion.

\vskip .2 cm

\paragraph{\bf $k$ preserves the orientation} We noted above that
$k$ maps the orbit $K\{ e\}$ onto $K\{ f\}$ without specifying an orientation,
and we noted that we may assume that $k(e)=f$
preserving orientation. It is a standard fact that $k:S^1\to
S^1$ extends to a differentiable self map $\tilde{k}$ of $\D$ which
sends complementary triangles of $\tau_{*}$ onto complementary
triangles of $\tau$ \cite{P1}. If $\alpha$ is a differentiable
path between $e$ and $e'$ as above, then $\tilde{k}(\alpha )$ is a
differentiable path between $f$ and $f'':=k(e')$ which satisfies
the required properties. Note that it is not necessarily true that
$f'=\gamma'(f)$ and $f''$ are equal. However, the inductive
definition of the characteristic map $k$ immediately implies that
$k:e'\mapsto f''$ is orientation preserving.

\vskip .2 cm

\paragraph{\bf $G$ is orientation preserving} Elements
of $G$ are covering transformations for the surface $\D /K$. We show
that $\gamma:e\mapsto e':=\gamma (e)$ is orientation preserving,
and a similar statement for $f$ follows immediately. Denote by
$\pi_K:\D\to\D /K$ the universal covering map. Let $\alpha'$ be a
differentiable curve on $\D /K$ representing the covering transformation $\gamma$ which is transverse to $\pi_K(e)$
and crosses the minimal number of edges of $\tau_{*}/K$.
We denote by
$\alpha$ a part of the lift of $\alpha'$ to $\D$ which connects
$e$ and $e'$, so $\gamma (e\cap\alpha )=e'\cap\alpha$. Since
$\gamma$ preserves the orientation of $\alpha$, it follows that
$\gamma :e\mapsto e'$ is orientation preserving.

\vskip .2 cm

\paragraph{\bf $k$ conjugates $K$ onto itself.}
Recall that $k$ conjugates a finite-index subgroup $H$ of
$PSL_2(\mathbb{Z})$ onto $K$. Since $G$ preserves the orientation
of the orbits $K\{ e\}$ and $K\{ f\}$ in the sense of the previous
paragraph and $k$ maps $K\{ e\}$ onto
$K\{ f\}$, it follows that $k$
conjugates the action of $K$ on the orbit $K\{ e\}$ onto the
action of $K$ on the orbit $K\{ f\}$. Since $H$ and $K$ have the
same index in $PSL_2(\mathbb{Z})$ it follows that $H=K$. $\Box$

\vskip .4 cm

Thus, $k$ descends to a self map $\bar{k}$ of $\D /K$ sending
the tesselation $\tau_{*}/K$ onto the tesselation $\tau /K$ and vice
versa. Let $\bar{e}\in\tau_{*}/K$ and $\bar{f}\in\tau /K$ be the
corresponding edges on $\D /K$ of the orbit $K\{
e\}\subset\tau_{*}$ and its corresponding orbit $K\{ f\}\in\tau
-\tau_{*}$ under the Whitehead move defining the edge $E$.

\vskip .2 cm

It follows from the proof above that $\bar{e}$ is necessarily
mapped onto $\bar{f}$ by $\bar{k}$, whence
$\bar{k}^2(\bar{e})=\bar{e}$ and $\bar{k}^2(\bar{f})=\bar{f}$ with
the orientations of $\bar{e}$ and $\bar{f}$ reversed. This implies
that $\bar{k}^4(\bar{e})=\bar{e}$ with an orientation of $\bar{e}$
preserved. Since in addition $\bar{k}^4(\tau_{*}/K)=\tau_{*}/K$,
we conclude that $\bar{k}^4=id$. This implies that $k$ (after possibly
pre-composing by an element of $K$ and for simplicity
renaming the composition again by $k$) maps $e$ onto $f$
sending $\tau_{*}$ onto $\tau$, and vice versa. Since $k^2\in
\Gamma (E)$ maps $e$ onto itself by reversing its orientation,
we conclude that $k^2\in K'$ is an involution with fixed point on
$e$. Thus,
$$
\Gamma (E)~=~<K',k>
$$
where $K'<PSL_2(\mathbb{Z})$, and $k\in Mod(\S )$ with $k^2\in
K'-K$ and $k^4=id$. In particular, $k^2$ is an elliptic involution
whose fixed point lies on $e$ and $K'\neq K$ if $k^2$ is non trivial. Note
that any Whitehead move on the once punctured torus can be
obtained as a homeomorphism of the torus which interchanges the
two tesselations, and hence the corresponding edge is inverted. In the
following example we show that edges admit orientation reversing isotropy
also for higher
genus.

\vskip .2 cm

\paragraph{\bf Example 4.2} We give in Figure~1 just one illustrative
example of  a surface $\D /K$ with a
distinguished quadrilateral, together with a self homeomorphism
performing a Whitehead move on the quadrilateral. In this figure,
the homeomorphism $h$ is a rotation by $\pi/4$ along the
horizontal axis. The dots represent the punctures of the surface.
\begin{figure}[htbp]
\centering
\includegraphics[scale=0.4]{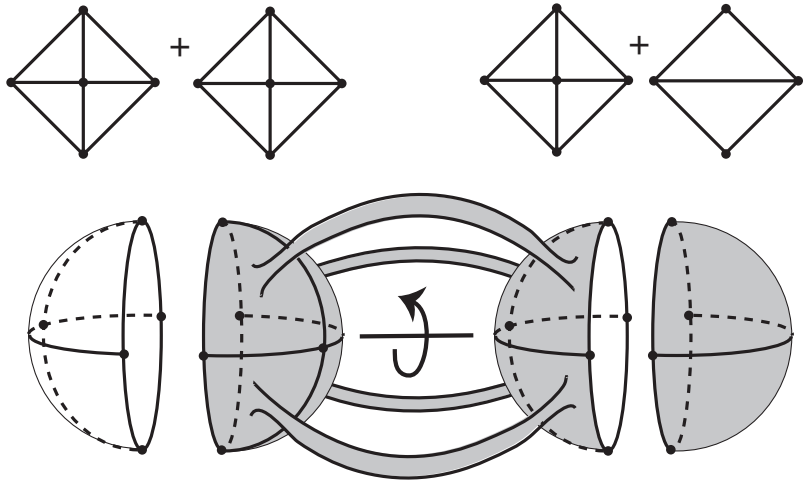}
\caption{Inverting an edge} \label{Rotation}
\end{figure}

\vskip .2 cm

We also note that there are infinitely many edges $E\in\mathcal
X$ with $\Gamma (E)-\Gamma ^{+}(E)=\emptyset$. This follows from the fact
that there are infinitely many Whitehead moves on finite
surfaces (whose Euler characteristics are increasing without
bound) such that there is no homeomorphism of the surface which
maps the starting tesselation onto the ending tesselation and by
Lemma 4.1.

\vskip .2 cm

Consider a two-cell of $\mathcal X$ with one vertex at the
basepoint $\tau_{*}$. Recall that for pentagon and square two-cells, each vertex is invariant under a finite-index subgroup $K$
of $PSL_2(\mathbb{Z})$, and edges correspond to Whitehead moves
invariant under this group $K$. For coset two-cells, either the long edge
has  $\tau_{*}$ as endpoint,  in which case all edges are Whitehead
moves invariant under $K$, or the long edge does have
$\tau_{*}$ as endpoint, in which case the long edge is given by a generalized
Whitehead move.

\vskip .2 cm

We may now apply Brown's theorem \cite[Theorem 1]{Bro} to obtain a
presentation of the modular group $Mod(\S )$ since  it  acts cellularly on the connected and simply connected
triangulation complex $\mathcal{X}$ with a single vertex orbit.
In fact, we shall introduce a
somewhat larger set of generators than necessary for the
application of Brown's theorem in order to obtain a simpler
presentation.

\vskip .2cm

It is a standard fact (which follows from Tietze's Theorem for instance) that for a given presentation,
if one adds extra generators, then an equivalent presentation arises by expressing the new generators
in terms of the old as new relations. One can in effect replace any
occurrence of a subsequence of old generators in the old relations
by new generators, in order to presumably simplify the presentation.
We shall ultimately give the presentation of an abstract group  ${\mathcal{G}}$, which is
equivalent in this sense to the presentation of the group in Brown's Theorem.

\vskip .2cm

The
set of edges of $\mathcal X$ that are not inverted by the action
of $Mod(\S )$ can be oriented consistently for the action of
$Mod(\S )$, and we fix one such orientation on each such edge.
Let $\mathcal{E}^{+}$
be the set of edges which are not inverted by the action of $Mod(\S
)$ that have initial point $\tau_{*}$. If $E=(\tau_{*},\tau
)\in\mathcal{E}^{+}$ then $\tau$ is invariant under a finite-index
subgroup $K$ of $PSL_2(\mathbb{Z})$.
By Lemma 4.1 and the subsequent
discussion, there is no homeomorphism $f:\D /K\to\D /K$ such that
$f(\tau_{*}/K)=\tau /K$ and $f(\tau /K)=\tau_{*} /K$.
The elements of $\mathcal{E}^{+}$ therefore are obtained by taking all
finite-index torsion free subgroups of $PSL_2(\mathbb{Z})$ and
performing all possible Whitehead moves on $\tau_{*}$ invariant
under the chosen groups, where the Farey tesselation $\tau_{*}$
and the image tesselation satisfy the additional property of not being
mapped onto each other by a single map conjugating the group onto
itself.
The images of $\tau_{*}$ under the Whitehead moves are the
terminal vertices of edges in $\mathcal{E}^{+}$. For any such
$E=(\tau_{*},\tau )\in\mathcal{E}^{+}$, we fix the characteristic
map $g_{E}\in Mod(\S )$ such that $g_E(\tau_{*})=\tau$ and the
standard distinguished oriented edge $e_0=(-1,1)$ of $\tau_{*}$ is
mapped to either itself if the Whitehead move is not along an
orbit of $e_0$, or it is mapped onto $f_0=(-i,i)$ if the Whitehead
move is along an orbit of $e_0$. The characteristic map is
uniquely determined by these conditions, and we fix this choice
$g_E$. (Notice that $\mathcal{E}^{+}$ is larger than necessary, since
it is enough to take only the edges corresponding to
representatives of {\sl conjugacy classes} in $PSL_2(\mathbb{Z})$ of
finite-index subgroups. However, this larger set simplifies the presentation,
and not much is lost because both sets are infinite.)

\vskip .2 cm

Let $\mathcal{E}^{-}$ denote the set of inverted edges with initial point
$\tau_{*}$.  By Lemma 4.1, an edge $E=(\tau_{*},\tau )$ is
inverted if there exists $k:S^1\to S^1$ such that
$k(\tau_{*})=\tau$, $k(\tau )=\tau_{*}$ and $kHk^{-1}=H$, where
$\tau$ is invariant under a torsion free finite-index subgroup $H$
of $PSL_2(\mathbb{Z})$. The isotropy group $\Gamma (E)$ of the
cell underlying $E\in{\mathcal E}^{-}$ is the subgroup of $Mod(\S )$
generated by $H'$ and $k$, where $H'>H$ is the maximal subgroup of
$PSL_2(\mathbb{Z})$ under which $\tau$ is invariant, where $k^2\in
PSL_2(\mathbb{Z})$ and where $k^4=id$. Fix some choice
$g_E$ of characteristic map associated to $E$ and take
$k=g_E$. (Again, we take $\mathcal{E}^{-}$ larger then necessary
for ease in writing down the relations.)

\vskip .2 cm

Denote by $\mathcal F$ the set of two-cells of $\mathcal X$
based at $\tau_{*}$. This condition implies that for each coset
two-cell in $\mathcal F$ the initial point of the long edge is
$\tau_{*}$ and each edge of the cell is consequently obtained by
a geometric Whitehead move, i.e., one invariant under a subgroup of
$PSL_2(\mathbb{Z})$ as opposed to a conjugate of a subgroup of
$PSL_2(\mathbb{Z})$.
This property also holds for pentagon and
square cells in $\mathcal{F}$. The set $\mathcal F$ is obtained by
taking all torsion free finite-index subgroups of
$PSL_2(\mathbb{Z})$ and taking all possible pentagon, square and
coset edge-paths in $\mathcal X$ corresponding to the chosen
groups.

\vskip .2cm

Note that a single choice of group for coset two-cells
gives countably many coset cells because there are countably
many finite-index subgroups and each finite-index subgroup yields
finitely many coset cells. In the case of a square or a pentagon
cell, a choice of a finite-index group determines finitely many
cells because there are finitely many edge orbits in $\tau_{*}$
under the group. (Again, we could have taken only representatives
of orbits of two-cells based at $\tau_{*}$ together with subsets
of $\mathcal{E}^{\pm}$, but for
simplicity later, we have expanded these sets.)

\vskip .2 cm

Let us assume for a moment that we had instead chosen for
$\mathcal{E}^{\pm}$ and $\mathcal F$ only representatives of
classes under the action of $Mod(\S )$. We describe the assignment
of a composition of elements in $Mod(\S )$ (depending upon the above
choices) to any closed boundary edge-path of a two-cell in
$\mathcal F$ in order to give a relation corresponding to the two-cell as in Brown's theorem \cite{Bro}.
Given an oriented edge
$E=(\tau_{*},\tau )$ starting at $\tau_{*}$, we assign to it a chosen
element of $g\in Mod(\S )$ such that $\tau =g(\tau_{*})$. If
$E\in\mathcal{E}^{+}$ then set $g:=g_E$. If
$E\notin\mathcal{E}^{+}$ is not inverted by $Mod(\S )$ then
$g:=\gamma\circ g_{E'}$, where $E'=(\tau_{*},\tau
')\in\mathcal{E}^{+}$, and $\gamma\in PSL_2(\mathbb{Z})$ satisfies
$\gamma (E')=E$; $\gamma$ is well-defined up to
pre-composition by an element of $\Gamma (E')$. The two Whitehead
moves from $\tau_{*}$ to $\tau$ and $\tau '$ determine
distinguished oriented edges $e$ and $e'$ of $\tau$ and $\tau '$,
and we choose unique $\gamma\in PSL_2(\mathbb{Z})$ mapping
$e'$ to $e$. Our choice of $\gamma\circ g_E$ is in this case
unique. If $E\in\mathcal{E}^{-}$, then set $g:=g_E\in \Gamma (E)-\Gamma ^{+}(E)$. If $E\notin
\mathcal{E}^{-}$ is inverted by the action of $Mod(\S )$ then
$g:=\gamma\circ g_{E'}$, where $E'\in\mathcal{E}^{-}$ and
$\gamma\in PSL_2(\mathbb{Z})$ with $\gamma (\tau ')=\tau$ and
$\gamma (e')=e$. The edge $E$ therefore ends at $g(\tau_{*})$, but
 it seems complicated to explicitly determine $E'$ and
$\gamma$.

\vskip .2 cm

Continuing to assume that we had chosen for
$\mathcal{E}^{\pm}$ and $\mathcal F$ only representatives of
classes under the action of $Mod(\S )$, consider a closed path
$\alpha$ of edges in $\mathcal X$ based at $\tau_{*}$.
Let
$(E_1,E_2,\ldots ,E_n)$ be the sequential edges of $\alpha$.
Denote by $g_1$ the unique element of $Mod(\S )$ chosen for the
edge $E_1=(\tau_{*},\tau_1)$ starting at $\tau_{*}$ as above, so
$g_1(\tau_{*})=\tau_1$. The edge $E_2=(\tau_1 ,\tau_2)$
is therefore of the form $g_1(E_2')$ for an edge $E_2'=(\tau_{*},\tau_2')$
based at $\tau_{*}$.
Denote by $g_2$ the unique element of $Mod(\S
)$ associated to $E_2'$ as above, so $g_1\circ
g_2(\tau_{*})=\tau_2$. This implies that $E_3=(\tau_2,\tau_3)$ is
given by $g_1\circ g_2(E_3')$, where $E_3'$ starts at $\tau_{*}$.
Take $g_3\in Mod(\S )$ associated to the edge $E_3'$, and continue
in this manner until we exhaust all edges of $\alpha$. This yields
a composition $g_1\circ g_2\circ\ldots \circ g_n$ in terms of
generators such that $g_1\circ g_2\circ\ldots \circ
g_n(\tau_{*})=\tau_{*}$.  Thus, there is $\gamma\in
PSL_2(\mathbb{Z})$ such that $g_1\circ g_2\circ\ldots \circ
g_n=\gamma$, and this is the relation associated with a closed edge
path $\alpha$ based at $\tau_{*}$. It seems complicated to
determine the maps $g_i$ from the given description or to decide
which elements $\gamma\in PSL_2(\mathbb{Z})$ arise.
However, the choice of $g_i$ simplifies if we allow all edges with
initial point $\tau_{*}$, and this will require additional
relations as discussed before.

\vskip .2 cm

>From this point on, we go back to our choice of
$\mathcal{E}^{\pm}$ to consist of all edges with initial point
$\tau_{*}$ and of $\mathcal{F}$ to consists of all two-cells based
at $\tau_{*}$.

\vskip .2cm

We describe the relations associated to boundaries of
two-cells in $\mathcal{F}$. Let us start with a pentagon two-cell
$P$ based at $\tau_{*}$ whose boundary edges are
$\{E_1=(\tau_{*},\tau_1),E_2=(\tau_1,\tau_2),\ldots
,E_5=(\tau_4,\tau_{*})\}$. The pentagon two-cell $P$ is given by
changing an orbit of two adjacent edges $e_1,e_2$ of $\tau_{*}$
under a torsion free finite-index subgroup $K$ of
$PSL_2(\mathbb{Z})$ of index at least $9$. Assume first that
the distinguished oriented edge $e_0=(-1,1)$ of $\tau_{*}$ is not
an element of the orbit $K\{e_1,e_2\}$ and apply the algorithm of
Brown to get the edge-path relation, but using our extended set of
generators to simplify it.
We denote by $g_i$ the element of
$Mod(\S )$ which corresponds to the edge $E_i$. The first edge
$E_1$ gives $g_1:=g_{E_1}$, so $g_1(e_0)=e_0$. We find
$g_2:=\gamma\circ g_{E_2''}$, where $E_2''$ is a
representative of the orbit of $(\tau_{*},g_1^{-1}(\tau_2))$ and
$\gamma (E_2'')=(\tau_{*},g_1^{-1}(\tau_2))$ is chosen from
$PSL_2(\mathbb{Z})$ in a unique way as above (i.e., $\gamma
(e_0)=e_0$). However, since $\mathcal{E}^{+}$ consists of all edge
starting at $\tau_{*}$ we immediately obtain that $g_2:=g_{E_2'}$,
where $E_2':=(\tau_{*},g_1^{-1}(\tau_2))$. We likewise
obtain $g_i:=g_{E_i'}$, for $i=3,4,5$, where $E_3':=(\tau_{*},(g_1\circ
g_2)^{-1}(\tau_3))$,
$E_4':=(\tau_{*},(g_1\circ g_2\circ g_3)^{-1}(\tau_4))$, and
$E_5':=(\tau_{*},(g_1\circ g_2\circ g_3\circ
g_4)^{-1}(\tau_{*}))$. Under our assumption that $e_0\notin
K\{e_1,e_2\}$, we find $g_i(e_0)=e_0$ for $i=1,2,\ldots ,5$.
The relation associated  to $P$ is therefore
\begin{equation}
\label{pentagon} g_1\circ g_2\circ\cdots\circ g_5=id.
\end{equation}
On the other hand, now assume $e_0\in K\{ e_1\}$ and without loss
of generality we can assume that $e_1=e_0$. We choose $g_i$ as
above and note that $g_1:(\tau_{*},e_0)\mapsto (\tau_1,e_0')$,
where $e_0'=(-i,i)$ is the image of $e_0=(-1,1)$ under the
Whitehead move corresponding to $E_1=(\tau_{*},\tau_1)$, $g_1\circ
g_2:(\tau_{*},e_0)\mapsto (\tau_2,e_0')$, $g_1\circ g_2\circ
g_3:(\tau_{*},e_0)\mapsto (\tau_3,e_0'')$ where $e_0''$ is the
image of $e_0'$ under the Whitehead move corresponding to $E_3$,
$g_1\circ g_2\circ g_3\circ g_4:(\tau_{*},e_0)\mapsto
(\tau_4,e_0'')$ and $g_1\circ g_2\circ\cdots\circ
g_5:(\tau_{*},e_0)\mapsto (\tau_{*},\bar{e}_2)$ where $\bar{e}_2$
is the oriented edge $e_2$ with orientation given such that the
terminal point of $e_0$ is the initial point of $\bar{e}_2$.
Denote by $\gamma_{e_0,\bar{e}_2}\in PSL_2(\mathbb{Z})$ the unique
element which maps $e_0$ onto $\bar{e}_2$ with the given
orientation. Thus,  $\gamma_{e_0,\bar{e}_2}$ is the composition of
the primitive parabolic element with fixed point at the terminal
point of $e_0$ which maps $e_0$ onto $e_2$ and the involution
which reverses $e_2$. We obtain the following relation
\begin{equation}
\label{pentagon1} g_1\circ g_2\circ\cdots\circ
g_5=\gamma_{e_0,\bar{e}_2}.
\end{equation}
When $e_2=e_0$ the relation is similarly
\begin{equation}
\label{pentagon2} g_1\circ g_2\circ\cdots\circ
g_5=\gamma_{e_0,\bar{e}_1}.
\end{equation}

\vskip .2 cm

Let $P$ be a square cell in $\mathcal F$. Assume that $P$ is
obtained by Whitehead moves along the nonadjacent orbits $K\{
e_1\}$ and $K\{e_2\}$ of edges $e_1,e_2$ in $\tau_{*}$, where $K$
is a torsion free finite-index subgroup of $PSL_2(\mathbb{Z})$. If
$e_0\notin K\{e_1,e_2\}$ then
\begin{equation}
\label{square} g_1\circ\cdots\circ g_4=id,
\end{equation}
where $g_i$ are chosen as above. If $e_i=e_0$ then we obtain a
relation
\begin{equation}
\label{square1} g_1\circ\cdots\circ g_4=s_{e_0},
\end{equation}
where $s_{e_0}\in PSL_2(\mathbb{Z})$ is the involution which
reverses $e_0$. The proofs of both relations for the square cell
$P$ depend upon keeping track of where  $e_0$ is mapped, and it is sufficiently similar
to the pentagon two-cell that we do not repeat it.

\vskip .2 cm

Let $P\in\mathcal F$ be a coset two-cell for the edge
$e\in\tau_{*}$ and for the groups $K_1<K<PSL_2(\mathbb{Z})$. If
$e\neq e_0$ then we obtain a relation
\begin{equation}
\label{coset} g_1\circ\cdots\circ g_k=id,
\end{equation}
where $k=[K:K_1]$ and $g_i$ are uniquely chosen as above. Note
that a single choice of $K_1<K$ gives a decomposition of the orbit
$K\{ e\}$ into $k$ disjoint coset orbits $K_1\{e_1\},K_1\{
e_2\},\ldots ,K_1\{ e_k\}$, where $e_i\in K\{ e\}$. This gives
$k!$ possible permutations on $K_1\{e_1\},K_1\{ e_2\},\ldots
,K_1\{ e_k\}$ which in turn produce $k!$ coset two-cells with
the long edge given by the Whitehead move on $K\{ e\}$. Note that
$g_1=g_E$ where $E=(\tau_{*},\tau )$ and $\tau$ is the image of
the Whitehead move along $K\{ e\}$. The other $g_i$, for
$i=2,3,\ldots ,k$, are given by the translation to $\tau_{*}$ of
the short edges. If $e=e_0$ then we obtain a relation
\begin{equation}
\label{coset1} g_1\circ\cdots\circ g_k=s_{e_0}.
\end{equation}

\vskip .2 cm

The desired group ${\mathcal{G}}$ is by definition the free product of
$PSL_2(\mathbb{Z})$, $\Gamma (E)=\Gamma ^{+}(E)$ for
$E\in\mathcal{E}^{+}$, $\Gamma (E)$ for $E\in  \mathcal {E }^{-}$ and a free
group generated by $g_E$ for $E\in\mathcal{E}^{+}$. The modular
group $Mod(\S )$ is the quotient of ${\mathcal{G}}$ by a set of
relations as follows.

\vskip .3 cm

\paragraph{\bf Theorem 4.3} {\it The modular group $Mod(\S )$ is
generated by the isotropy subgroup $PSL_2(\mathbb{Z})$ of the
basepoint $\tau_{*}\in\mathcal X$, the isotropy subgroups
$\Gamma (E)$ for $E\in\mathcal{E}^{\pm}$,
and by the elements $g_E$ for
$E\in\mathcal{E}^{+}$. The following relations on these
generators give a complete presentation of $Mod(\S )$:

\vskip .1 cm

\noindent a) The inclusions of $\Gamma (E)$ into $PSL_2(\mathbb{Z})$,
for $E\in\mathcal{E}^{+}$, are given by $\Gamma (E)=K'$,
where the terminal endpoint of $E$ is invariant under
the finite-index subgroup $K'<PSL_2(\mathbb{Z})$;

\vskip .1 cm

\noindent b) The inclusions of $\Gamma ^+(E)$ into $PSL_2(\mathbb{Z})$,
for $E\in\mathcal{E}^{-}$, are given by $\Gamma (E)=K'$,
where the terminal endpoint of $E$ is invariant under
the finite-index subgroup $K'<PSL_2(\mathbb{Z})$;

\vskip .1 cm

\noindent c) The relations introduced by the boundary edge-paths
of two-cells in $\mathcal F$ given by the equations}
(\ref{pentagon}), (\ref{pentagon1}), (\ref{pentagon2}),
(\ref{square}), (\ref{square1}), (\ref{coset}) {\it and}
(\ref{coset1});

\vskip .1 cm

{\it \noindent d) The redundancy relations: for any two edges $E$
and $E'$ in $\mathcal{E}^{\pm}$ and for any $\gamma\in
PSL_2(\mathbb{Z})$ such that $\gamma (E)=E'$, we get the relation
$$g_{E'}\circ \gamma '=\gamma\circ g_E,$$
 where $\gamma '$ is the unique
element of $PSL_2(\mathbb{Z})$ that satisfies $\gamma '(e_0)=e_1'$
with $e_1'=g_{E'}^{-1}( \gamma (e_0))$.}

\vskip .2 cm

\paragraph{\bf Proof} The fact about the generators of $Mod(\S )$
follows directly from Brown's theorem \cite{Bro} and from our
choice of $\mathcal{E}^{\pm}$ even larger than necessary. The
relations from \cite[Theorem 1]{Bro} are included in our theorem
as follows. The relations (i) are empty in our case. The relations
(ii), (iii) and (iv) translate easily to relations a), b) and c)
in our theorem, respectively. The relations d) are extra relations
needed because we have taken a larger set of generators than in
Brown's presentation. If $g_E(e_0)=e_0$, then the relation d) is
immediate. If $g_E(e_0)\neq e_0$, then $g_{E'}(e_0)=e_0$ (since
$\gamma\notin K$), and d) follows by reversing the roles of $g_E$
and $g_{E'}$. $\Box$

\section{No central elements}

In analogy to the case of surfaces of finite type, we have:

\vskip .2 cm

\paragraph{\bf Theorem 5.1} {\it The modular group $Mod(\S )$ of
the punctured solenoid $\S$ has trivial center.}

\vskip .2 cm

\paragraph{\bf Proof} Let $h\in Mod(\S )$ be a central element, so $h$
is a word in the generators of Theorem~4.3.
If $h\in
PSL_2(\mathbb{Z})<Mod(\S)$, then since
$PSL_2(\mathbb{Z})$ has trivial center, it
follows that $h=id$.

\vskip .2 cm

Assuming that $h\notin PSL_2(\mathbb{Z})$, we have
$h(\tau_{*})=\tau\neq\tau_{*}$.  Since $h$ is supposed to be central, we
have $g\circ h\circ g^{-1}=h$ for all $g\in Mod(\S )$. Taking $g\in PSL_2(\mathbb{Z})$,
we find $g\circ h
(\tau_{*})=h(\tau_{*})$, i.e., $g(\tau )=\tau$ for all
$g\in PSL_2(\mathbb{Z})$. By Proposition 5.2 below, we conclude $\tau
=\tau_{*}$, i.e., $h(\tau_{*})=\tau_{*}$, which implies
$h\in PSL_2(\mathbb{Z})$. This gives a contradiction, and so again $h=id$.
$\Box$

\vskip .2 cm

\paragraph{\bf Proposition 5.2} {\it A TLC
tesselation of $\D$ is invariant under $PSL_2(\mathbb{Z})$ if and
only if it is the Farey tesselation $\tau_{*}$.}

\vskip .2 cm

\paragraph{\bf Proof} The implication that $\tau _*$ is the unique
tesselation invariant under $PSL_2({\mathbb Z})$ is given in \cite[Lemma 7.3]{PS}.  For
the converse,
recall that $PSL_2(\mathbb{Z})$ contains order
two elliptic elements with fixed points on each edge of $\tau_{*}$
and order three elliptic elements with fixed points at the center
of each ideal complementary triangle of $\tau_{*}$. Let $T$ be a
complementary triangle of $\tau$ containing the fixed point $a$
of an elliptic element $\gamma\in PSL_2(\mathbb{Z})$ of order
three. It is an exercise in elementary hyperbolic geometry to show
that $a$ is the center of $T$.

\vskip .2 cm

Let $b$ be a fixed point of an elliptic involution $\gamma\in
PSL_2(\mathbb{Z})$. If $b$ is in the interior of a complementary
triangle $T$ of $\tau$ then $\gamma (T)\neq T$ and $\gamma (T)\cap
T\neq\emptyset$. Thus, the image of the boundary of $T$ under $\gamma$
intersects transversely the boundary of $T$. This is in
contradiction to the assumption that $\gamma$ fixes $\tau$.
It follows that $b$ must lie on an edge of $\tau$.

\vskip .2 cm

Let $a$ be the fixed point of an elliptic element of
$PSL_2(\mathbb{Z})$ of order three and let $b_1,b_2,b_3$ be fixed
points of three elliptic involutions of $PSL_2(\mathbb{Z})$ that
are shortest distance to $a$ among all such involutions. Thus,
$b_1,b_2,b_3$ lie on a hyperbolic circle centered at $a$. If
$T$ is the ideal triangle in the complement of $\tau$ whose center
is $a$, then the boundary sides of $T$ are tangent to this
circle. Since $b_1,b_2,b_3$ must lie on edges of $\tau$, this
implies that the boundary sides of $T$ are tangent at the points
$b_1,b_2,b_3$. It follows that $T$ is a
complementary triangle of $\tau_{*}$ as well. Since this is true
for an arbitrary $T$, it follows that indeed $\tau =\tau_{*}$. $\Box$

\vskip .4 cm

We consider the action of $Mod(\S )$ on the first barycentric
subdivision $\mathcal{X}'$ of $\mathcal{X}$.

\vskip .2 cm

\paragraph{\bf Proposition 5.3} {\it The first barycentric subdivision ${\mathcal X}'$
of $\mathcal X$ is a simplicial complex on which $Mod(\S)$ acts
simplicially.}

\vskip .2 cm

\paragraph{\bf Proof} Note that $Mod(\S )$ preserves cells of $\mathcal X$. An isotropy
group of a vertex of $\mathcal X$ is a conjugate of
$PSL_2(\mathbb{Z})$. We showed that the isotropy group of an edge
is either a finite-index subgroup of $PSL_2(\mathbb{Z})$ which
preserves the orientation of the edge or it is generated by an
element of $Mod(\S )$ which reverses the orientation of the edge
and by a finite-index subgroup of $PSL_2(\mathbb{Z})$ which
preserves the orientation of the edge. In the first case, each
element of the isotropy group fixes each point on the edge. In the
second case, an element either fixes each point of the edge or
fixes the midpoint and reflects the endpoints of the edge.

\vskip .2 cm

Let $C$ be a coset two-cell with long edge given by a Whitehead
move on TLC tesselation $\tau_{*}$ invariant under
$K<PSL_2(\mathbb{Z})$ along the orbit of $e\in\tau$ and with the
short edges given by Whitehead moves invariant under a
subgroup $H<K$. The isotropy subgroup of $C$ is a finite
extension in $PSL_2(\mathbb{Z})$ of $H$. To prove the claim,
it is enough to show that the long edge cannot be mapped onto a
short edge, and this is true because the rectangles in which change of
diagonals for the Whitehead move occur must be mapped onto the
rectangles on which change of diagonals occur.  However, the two
groups have different indexes in $PSL_2(\mathbb{Z})$ which gives a
contradiction, and the isotropy group $\Gamma (C)$ therefore acts by
fixing each point in $C$. Since an arbitrary coset two-cell is the
image of some $C$ as above, the same statement holds for an
arbitrary coset two-cell.

\vskip .2 cm

Let $P$ be a pentagon two-cell based at $\tau_{*}$ obtained by
Whitehead moves along $K\{ e_1,e_2\}$. Thus, $P$ has a subgroup
$K'$, where $PSL_2(\mathbb{Z})>K'>K$, of its isotropy group $\Gamma (P)$
fixing each point of $P$. If $\Gamma (P)\neq K'$, then it is
generated by $K'$ and a single element of $h\in Mod(\S )$ which maps the first
edge onto the second. Since $h^5(\tau_{*})=\tau_{*}$,
we conclude $h^5\in PSL_2(\mathbb{Z})$. Thus, $h^5$ is mapping
class like and therefore $h$ is mapping class like. The map $h$
fixes the center of $P$ and rotates by the angle $2\pi /5$ the
pentagon $P$. The situation for a pentagon not based at $\tau_{*}$
is the same.

\vskip .2 cm

Finally, suppose $Q$ is a square two-cell obtained by Whitehead moves
along $K\{ e_1,e_2\}$ then $\Gamma (Q)>K'$, where
$K<K'<PSL_2(\mathbb{Z})$. It is possible  a priori that
$\Gamma (Q)\neq K'$, in which case the elements $h\in \Gamma (Q)-K'$
permute edges of $Q$ and fix the center of $Q$.~~~$\Box$

\vskip .4 cm

We finally investigate the topological fundamental group
$\pi_1(\mathcal{Y})$ of the quotient space $\mathcal Y
=\mathcal{X}/Mod(\S )$. To begin, we describe a natural surjection
$\phi :Mod(\S )\to \pi_1(\mathcal{Y})$ as follows. Denote by $\Pi
:\mathcal{X}\to\mathcal{Y}$ the quotient map, let $h\in Mod(\S )$
be arbitrary, and define $\tau =h(\tau_{*})$. Let $\gamma$ be an
edge-path between $\tau_{*}$ and $\tau$ in $\mathcal{X}$. and
define
$$
\phi (h):=[\Pi (\gamma )],
$$
where $[\Pi (\gamma )]$ is the homotopy class of the closed curve
$\Pi (\gamma )$ based at $\Pi (\tau_{*})$, i.e., $[\Pi (\gamma
)]\in \pi_1(\mathcal{Y},\Pi (\tau_{*}))=\pi_1(\mathcal{Y})$. It is
a standard fact that $\phi$ is a well-defined and surjective
homomorphism.

\vskip .2cm

Let
$\mathcal N$ be the group generated by the isotropy subgroups of all
vertices of $\mathcal X'$, where $\mathcal X'$ denotes the first barycentric
subdivision of $\mathcal X$, so $\mathcal N$ is normal
in $Mod(\S )$. In fact, $\mathcal N$ is generated by all
conjugates of $PSL_2(\mathbb{Z})$ and by the isotropy groups of
edges and  two-cells of $\mathcal X$.
The isotropy group of an edge in $\mathcal X$ fixes the
center of the edge and therefore belongs to the isotropy group of
a vertex in $\mathcal X'$, and likewise the isotropy group of a two-cell in $\mathcal X$ fixes a vertex of $\mathcal X'$. Moreover,
any element of $Mod(\S )$ which fixes a point in $\mathcal X'$
fixes a point in $\mathcal X$. By Proposition 5.3 and by a standard
result \cite{Arm} we get:

\vskip .2 cm

\paragraph{\bf Theorem 5.4} {\it The topological fundamental group of
${\mathcal Y}={\mathcal X}/Mod(\S )$ satisfies
$$
\pi_1({\mathcal Y})= Mod(\S )/\mathcal{N},
$$
where $\mathcal{N}<Mod(\S )$ is generated by the isotropy groups
of vertices, edges and two-cells of $\mathcal{X}$.} $\Box$

\vskip .2 cm

By our discussion above, each element of $Mod(\S )$ which fixes a
cell in $\mathcal{X}$ is mapping class like, i.e. it conjugates a
finite index subgroup of $G$ onto itself. Therefore $\mathcal{N}$
is generated by some mapping class like elements. We pose the
following question:

\vskip .2 cm

\paragraph{\bf Question} Is $\mathcal{N}$ equal to the normal
subgroup of $Mod(\S )$ generated by all mapping class like
elements?

\end{document}

%% file: imsmark.tex
\def\IMSmarkvadjust{0 pt}
\def\IMSmarkhadjust{0 pt}
\def\IMSmarkhpadding{0 pt}
\def\IMSpubltext{Published in modified form:}
\def\SBIMSMark#1#2#3{
 \font\SBF=cmss10 at 10 true pt
 \font\SBI=cmssi10 at 10 true pt
 \setbox0=\hbox{\SBF \hbox to \IMSmarkhpadding{\relax}
                Stony Brook IMS Preprint \##1}
 \setbox2=\hbox to \wd0{\hfil \SBI #2}
 \setbox4=\hbox to \wd0{\hfil \SBI #3}
 \setbox6=\hbox to \wd0{\hss
             \vbox{\hsize=\wd0 \parskip=0pt \baselineskip=10 true pt
                   \copy0 \break%
                   \copy2 \break%
                   \copy4 \break}}
 \dimen0=\ht6   \advance\dimen0 by \vsize \advance\dimen0 by 8 true pt
                \advance\dimen0 by -\pagetotal
	        \advance\dimen0 by \IMSmarkvadjust
 \dimen2=\hsize \advance\dimen2 by .25 true in
	        \advance\dimen2 by \IMSmarkhadjust

%
%
  \openin2=publishd.tex
  \ifeof2\setbox0=\hbox to 0pt{}
  \else 
     \setbox0=\hbox to 3.1 true in{
                \vbox to \ht6{\hsize=3 true in \parskip=0pt  \noindent  
                {\SBI \IMSpubltext}\hfil\break
                \input publishd.tex 
                \vfill}}
  \fi
  \closein2
  \ht0=0pt \dp0=0pt
 \ht6=0pt \dp6=0pt
 \setbox8=\vbox to \dimen0{\vfill \hbox to \dimen2{\copy0 \hss \copy6}}
 \ht8=0pt \dp8=0pt \wd8=0pt
 \copy8
 \message{*** Stony Brook IMS Preprint #1, #2. #3 ***}
}